\documentclass[a4paper,12pt]{article}

\usepackage{amsmath}
\usepackage{amsfonts}
\usepackage{amssymb}
\usepackage{epsfig}

\newtheorem{thm}{Theorem}
\newtheorem{wn}[thm]{Corollary}
\newtheorem{lem}[thm]{Lemma}
\newtheorem{rem}[thm]{Remark}
\newtheorem{df}[thm]{Definition}
\newtheorem{pro}[thm]{Proposition}

\begin{document}
\def\thethm{\thesection.\arabic{thm}}
\renewcommand \theequation {\thesection.\arabic{equation}}
\title{The fixed point property and unbounded sets in
spaces of negative curvature}
\author{Bozena Piatek}
\date{}
\maketitle
\begin{abstract}
Motivated by the well-known cases of the real Hilbert ball and complete $\mathbb{R}$-trees, 
being both particular cases of CAT(-1) spaces, we give
an affirmative answer to the question of whether the geodesically boundedness property is 
a necessary and sufficient condition for a closed convex
subset K of a complete CAT($\kappa$) space, with $\kappa$ < 0, to have the fixed point
property for nonexpansive mappings
\end{abstract}

\section{Introduction}
\setcounter{equation}{0}\setcounter{thm}{0}

In metric fixed point theory we are searching for the conditions, the imposing of which on metric spaces or self-maps on them, implies the existence of fixed points for such maps.
Especially, we are concerned with contractions and nonexpansive mappings. Unlike, for example to the well-know Banach contraction principle, the solution of the problem in case of nonexpansive mappings is based on the geometry of spaces. If we assume that $X$ is a~Hilbert space and $K$ is a~closed and convex subset of $X$, then $K$ has a~fixed point property if $K$ is bounded -- this is the well known result due to G\"{o}hde. The analogous result for a~some special family of Banach spaces was proved by Browder and Kirk, independently, also in 1965. The generous exposition of this topic the reader may find in \cite{GK,GR}.
In eighties Ray improved the result of Browder et al. showing that the boundedness of $K$ is also a~necessary condition for the fixed point property (see \cite{ORay,Sine}). The analogous result with much shorter proof was proposed by Sine in \cite{Sine}. Very recently this result was generalized for the case of Banach spaces satisfying some additional geometrical conditions by  Takahashi, Yao and Kohsaka (see \cite{Kohsaka}). At the same time in \cite{Domingues} Dom\'inguez gave an~example of Banach space for which the result does not hold. 

A~quite different situation takes place if $X$ is assumed to be a~space of negative curvature. In \cite{GR} it was shown that any closed and convex subset $K$ of the real Hilbert ball $B$ (being an~infinite dimentional Hadamard manifold of constant curvature equal to ''$-1$'', see \cite{GR}, for curvature compare \cite[Section II.10]{BH}) has the fixed point property for nonexpansive mappings if and only if $K$ is geodesically bounded. The same result for complete $\mathbb{R}$--trees, which may be treated as spaces of constant curvature equal to ''$-\infty$'', was proved by Esp\'inola and Kirk in \cite{EK} (see also \cite{Kirk}).

Fixed points on CAT(0) spaces, or spaces of globally nonpositive curvature in the sense of Gromov, have been extensively studied in the last years by a~wide group of mathematicians (see for instance \cite{EF}, \cite{EP} or \cite{Kirk}). \cite{BH} provides a~very comprehensive exposition on CAT(0) spaces. Very recently some considerations on the fixed point property for unbounded subsets in case of complete CAT(0) spaces were shown by Esp\'inola and Piatek in \cite{EP}. In that paper among the others authors raised still an~open question on the fixed point property for unbounded subsets of complete CAT(0) spaces
with curvature bounded above and below by two negative numbers. Since the real Hilbert ball as well as $\mathbb{R}$--trees are special cases of CAT(-1) spaces our considerations lead to the more general question: whether the geodesically boundedness is necessary and sufficient condition for a~closed convex subset $K$ of a~complete CAT($\kappa$), $\kappa<0$ space to have a~fixed point property for nonexpansive mappings. The main goal of this paper is to give an~affirmative answer to this question. 

Our paper is organized in the following way. Section 2 provides some definitions and elementary properties of geometry of CAT(0) spaces. In Section 3 we formulate some technical lemmas which will be useful in the proof of our main theorem. The most interesting in this part seems to be Lemma 3.2 where we consider geometrical behavior of $\mathbb{H}^2$ showing another difference between CAT(0) and CAT($\kappa$) spaces with negative $\kappa$. The main result -- the fixed point property for geodesically bounded subsets of CAT($\kappa$) spaces with negative $\kappa$ is presented in the last section.

\section{Preliminaries}
\setcounter{equation}{0}\setcounter{thm}{0}

Let $(X,\rho)$ be a~geodesic metric space. $X$ is said to be {\it uniquely geodesic} if each pair of points is joining by a~unique metric segment which will be denoted by $[x,y]$ for $x,y\in X$. A~subset $K$ of $X$ is called {\it convex} if $[x,y]\subset K$ as long as $x,y\in K$ and {\it geodesically bounded} as long as there is no infinite geodesic in $K$.

Now we introduce the concept of model spaces $M_\kappa^2$, $\kappa \leq 0$ which we need to define CAT($\kappa$) spaces. In \cite{BH} the reader can find a~very generous exposition on CAT($\kappa$) spaces also in case if $\kappa$ is positive.

Let us consider the space $\mathbb{R}^3$ endowed with the symmetric bilinear form which associates to vectors
$u=(u_1,u_2,u_3)$ and $v=(v_1,v_2,v_3)$ the real number
$\langle u| v\rangle$ defined by
$$
\langle u|v\rangle = u_1v_1+u_2v_2-u_3v_3.
$$

Let $\mathbb{H}^2$ be a~set
$$
\{x=(x_1,x_2,x_3)\in\mathbb{R}^3 \ |\ \langle x,x\rangle=-1 \land x_3\geq 1 \}
$$
Then $\mathbb{H}^2$ with a~function $d\colon \mathbb{H}^2 \times \mathbb{H}^2 \to \mathbb{R}$ defined by
$$
d(u,v)=\mbox{arc}\cosh\langle u,v\rangle
$$
is a~geodesic space.
In a~similar way one may define $\mathbb{H}^n$ being a~subset of $\mathbb{R}^{n+1}$, $n\in\mathbb{N}$ and infinite dimentional $\mathbb{H}^\infty$ (a~subset of Hilbert space $l^2$) being isometric to the real Hilbert ball (for generous exposition of the real and complex Hilbert ball see \cite[Section II.32]{GR}).

The Model Spaces $M_{\kappa}^2$ for $\kappa \leq 0$ are defined as follows.

\begin{df}
Given $\kappa \in (-\infty,0]$, we denote by $M_\kappa^2$ the following metric
spaces:
\begin{itemize}
\item[(1)] if $\kappa =0$ then $M^2_0$ is the Euclidean space ${\mathbb
E}^2$;%
\item[(2)] if $\kappa <0$ then $M^2_{\kappa}$ is obtained from the hyperbolic
space ${\mathbb H}^2$ by multiplying the distance function by the
constant $1/\sqrt{-\kappa}$.
\end{itemize}
\end{df}

Let $X$ be a~geodesic space. {\it A~geodesic triangle}
$\triangle(x,y,z)$ is a~set consisting of three points
$x,y,z\in X$ and all metric segments $[x,y]$, $[y,z]$ and $[z,x]$ (the edges of $\triangle$). By {\it a~comparison triangle} on $M^2_\kappa$ we will understand the triangle
$\triangle(\bar{x},\bar{y},\bar{z})$ with the same lengths of edges.
The comparison triangle always exists and it is unique up to isometry.

A geodesic triangle $\triangle$ in $X$ is said to satisfy the
CAT$(\kappa)$ {\sl inequality} if, given $\bar{\triangle}$ a comparison
triangle on $M_{\kappa}^2$ for $\triangle$, for all $a,b\in \triangle$
$$
\rho(a,b)\leq d(\bar{a},\bar{b}),
$$
where $\bar{a},\bar{b}\in\bar{\triangle}$ are the {\it
comparison points} of $a,b$, respectively.

\begin{df}
$X$ is called a CAT$(\kappa)$ space, $\kappa\leq 0$, if $X$ is a
geodesic space such that all of its geodesic triangles satisfy the
CAT$(\kappa)$ inequality.
\end{df}

Next we present some properties of CAT($\kappa$) spaces (with $\kappa\leq 0$) but to do this we will need the notion of Alexandrov's angle. Let $(x,u,v)$ be a~triple in
$X$ and $(\bar{x},\bar{u},\bar{v})$ a~comparison triple on the Euclidean plane. Assume $u,v\neq x$. Then the {\it comparison angle} $\angle_x(u,v)\in[0,\pi]$ is the (Riemannian) angle at $\bar{x}$ subtended by the segments $[\bar{x},\bar{u}]$, $[\bar{x},\bar{v}]$. Now let $y,z$ be points in $X$ and let $\sigma\colon[0,d(x,y)]\to X$ and $\tau\colon[0,d(x,z)]\to X$ be the geodesics from $x$ to $y$ and $z$, respectively. Then the Alexandrov angle is defined as
$$
\angle_x(y,z):= \lim_{s^{\prime}, t^{\prime}\to 0}\angle_x(\sigma (s^{\prime}), \tau (t^{\prime})),
$$
if the previous limit exists.

\begin{pro}{\bf (compare \cite{BH})\\}\label{pro_main}
Let $X$ be a~CAT($\kappa$) space, $\kappa\leq 0$. Then:
\begin{enumerate}
 \item[\emph{(i)}] Each pair of points $x,y\in X$ is joining by a~unique metric segment $[x,y]$.
 \item[\emph{(ii)}] If $\Delta(\bar{x},\bar{y},\bar{z})$ is a~comparison triangle of $\Delta(x,y,z)$ on $M^2_\kappa$, then the Alexandrov angle $\angle_{x}(y,z)$ is well defined and satisfies
$$
\angle_{x}(x,y)\leq \angle_{\bar{x}}(\bar{y},\bar{z}).
$$
 \item[\emph{(iii)}] If $\kappa<0$, then $X$ is also a~CAT($\kappa^\prime$) space for all $\kappa< \kappa^\prime\leq 0$.
 \item[\emph{(iv)}] If $(x_n)$ is a~bounded subset of a~CAT(0) space $X$, then there is a~unique asymptotic center $A((x_n))$ satysfying
     $$
     \limsup_{n\to\infty}d(A((x_n)),x_n)=\min_{y\in X}\limsup_{n\to\infty}d(y,x_n).
     $$
 \item[\emph{(v)}] Each triangle $\Delta(\bar{x},\bar{y},\bar{z})$ in $\mathbb{H}^2$ satisfies so called hyperbolic cosine law, namely,
$$
\cosh d(\bar{y},\bar{z})=\cosh d(\bar{x},\bar{y})\cosh d(\bar{x},\bar{z})-\sinh d(\bar{x},\bar{y})\sinh d(\bar{x},\bar{z})\cos\angle_{\bar{x}}(\bar{y},\bar{z}),
$$
where by $\angle_{\bar{x}}(\bar{y},\bar{y})$ we understand the Alexandrov angle in $\mathbb{H}^2$.
 \item[\emph{(vi)}] If $X$ is a~CAT(-1) space, then for each triangle $\Delta(x,y,z)$ the following inequality
$$
\cosh \rho(y,z)\geq\cosh \rho(x,y)\cosh \rho(x,z)-\sinh \rho(x,y)\sinh \rho(x,z)\cos\angle_{x}(y,z)
$$
holds.
\end{enumerate}

\end{pro}

In the sequel we will also need the notion of the Busemann convexity. Namely, $X$ is said to be {\it Busemann convex} if for each pair of geodesics $\sigma\colon[0,l_1]\to X$ and $\tau\colon[0,l_2]\to X$ parametrized with respect to arc there is
$$
\rho(\sigma(tl_1),\tau(tl_2))\leq (1-t)\rho(\sigma(0),\tau(0)) + t \rho(\sigma(l_1),\tau(l_2)).
$$
Clearly, from the CAT(0) inequality and (iii) of Proposition \ref{pro_main} it follows that each CAT($\kappa$) space (with $\kappa \leq 0$) is the Busemann convex one.

${\mathbb R}$-trees are a particular class of CAT$(0)$ spaces with many applications in different fields. They are also referred to as spaces of ''$-\infty$'' constant curvature (see \cite[p. 167]{BH} for more details). The interested reader may check \cite{khamsi,EK,Kirk,mine} for recent advances on $\mathbb R$-trees and fixed points.

\begin{df}
An ${\mathbb R}$-tree is a~geodesic metric space $M$ such that:
 \begin{itemize}
 \item[(1)] for all $x,y\in M$ there is unique metric segment $[x,y]$ joining them;
 \item[(2)] if $x,y$ and $z\in M$ are such that $[y,x]\cap [x,z]=\{
 x\}$, then $[y,x]\cup [x,z]=[y,z]$.
 \end{itemize}
\end{df}

\section{Technical lemmas}
\setcounter{equation}{0}\setcounter{thm}{0}

In this section we propose some technical lemmas which will be useful in the proof of main theorems. We begin with a~basic behaviour of triangles on a~plane:
\begin{lem}
Let us fix $\alpha_0\in (0,2\pi)$ and consider a~triangle $\Delta(\bar{x},\bar{y},\bar{z})$ on the Euclidean plane with $d(\bar{x},\bar{y})=a+d$, $d(\bar{x},\bar{z})=b+d$, $d(\bar{y},\bar{z})=C$ and $\angle_{\bar{x}}(\bar{y},\bar{z})=\alpha\geq \alpha_0$. Let $\bar{u}$ and$ \bar{v}$ be chosen on $[\bar{x},\bar{y}]$ and $[\bar{x},\bar{z}]$, respectively, in such a~way that $d(\bar{u},\bar{y})=a$ ($d(\bar{v},\bar{y})=b$).

Hence
\begin{itemize}
\item[$\mbox{(i)}$] $C\to\infty$ if $a,b\to \infty$;
\item[$\mbox{(ii)}$] $
d(\bar{u},\bar{v}) < d(\bar{y},\bar{z})-d\sin^2\dfrac{\alpha_0}{2}
$ for $a,b$ large enough (with respect to $d$).
\end{itemize}
\end{lem}

%\begin{figure}[!htb]
%\mbox{\null}\hfill{\psfig{figure=najw1.eps,height=40mm}}\hfill\mbox{\null}
%\caption{Triangle for Lemma 1}
%\end{figure}

{\bf Proof.}

Let us denote $d(\bar{u},\bar{v})$ by $c$ and $h=C-c$. Then (i) is an~easy consequence of the cosine law on a~plane.
%Then the cosine law on a~plane implies that
%$$
%c^2=(a-b)^2+2ab(1-\cos\alpha)\geq(a-b)^2+2ab(1-\cos\alpha_0).
%$$
%At the same time
%$$
%C^2=c^2+2(a+b+d)d(1-\cos\alpha)\geq c^2+2(a+b+d)d(1-\cos\alpha_0),
%$$
%so $C>c$ and both tends to $\infty$ if $a,b\to \infty$.

Moreover,
$$
2ch+h^2=2(a+b+d)d(1-\cos\alpha),
$$
but $h\leq 2d$, so
$$
ch \geq (a+b)d(1-\cos\alpha)-d^2(1+\cos\alpha).
$$
On account of (i) and since $c\leq a+b$ we obtain
$$
h \geq d(1-\cos\alpha)-\dfrac{d^2}{c}(1+\cos\alpha) >  \dfrac{d}{2}\left(2\sin^2\dfrac{\alpha}{2}\right)\geq d\sin^2\dfrac{\alpha_0}{2}. \qquad\square
$$
\\

In next Lemma we will show how to estimate the length of third edge of a~triangle in a~CAT($-1$) space. This behaviour will be useful in Step 4 of the proof of our main Theorem. Let us notice that similar result in CAT(0) spaces does not hold.
\begin{lem}
Let $X$ be a~CAT(-1) space and consider a~sequence of triangles $(\Delta(x_n,y_n,z_n))_{n=1}^\infty$ such that
$$
\angle_{y_n}(x_n,z_n)\geq \dfrac{\pi}{2}, \qquad n\in\mathbb{N},
$$
and
$$
d(x_n,y_n)\to \infty, \qquad d(x_n,z_n)-d(x_n,y_n)\to 0 \qquad \mbox{for }n\to\infty.
$$
Then $d(y_n,z_n)\to 0$.
\end{lem}

%\begin{figure}[!htb]
%\mbox{\null}\hfill{\psfig{figure=najw2.eps,height=30mm}}\hfill\mbox{\null}
%\caption{Triangle for Lemma 2}
%\end{figure}

{\bf Proof.}
Le us consider a~triangle $\Delta(\bar{x},\bar{y},\bar{z})$ on $\mathbb{H}^2$ such that $d(\bar{x},\bar{y})=h$, $d(\bar{x},\bar{z})=h+\varepsilon$ and $\angle_{\bar{y}}(\bar{x},\bar{z})\geq \dfrac{\pi}{2}$. Then the hyperbolic cosine law implies
$$
\cosh d(\bar{x},\bar{y})\cosh d(\bar{y},\bar{z}) \leq \cosh d(\bar{x},\bar{z})
$$
from which it follows that
$$
\cosh d(\bar{y},\bar{z}) \leq \dfrac{\cosh (h+\varepsilon)}{\cosh (h)}\leq \dfrac{e^{h+\varepsilon}+e^{-\varepsilon-h}}{e^{h}}=e^{\varepsilon}+e^{-\varepsilon-2h}
$$
and finally
\begin{equation}\label{szaccosh}
d(\bar{y},\bar{z}) \leq \mbox{arc}\cosh\Big[e^{\varepsilon}+e^{-\varepsilon-2h}\Big].
\end{equation}

Now let us notice that the comparison triangle $\Delta(\bar{x}_n,\bar{y}_n,\bar{z}_n)$ of $\Delta(x_n,y_n,z_n)$ (with the same lenghts of edges and $\angle_{\bar{y}_n}(\bar{x}_n,\bar{z}_n) \geq \angle_{y_n}(x_n,z_n)$) satisfies \eqref{szaccosh}, but
$$
e^{d(x_n,z_n)-d(x_n,y_n)}+e^{-d(x_n,z_n)-d(x_n,y_n)} \to 1
$$
so
$$
d(y_n,z_n)\to 0.\qquad\square
$$
\\

\begin{rem}
As it was mentioned earlier if one considers comparison triangles on the Euclidean plane the previous Lemma is not true.
\end{rem}

\begin{pro}{\bf (see \cite[Corollary 5.5]{EP})\\}
Let $X$ be a~CAT(-1) space and let $x_0$, $x$ and $y \in X$ such that there exists $r$, $\varepsilon > 0$ with $d(u, v)\geq\varepsilon$,
where $u$ and $v$ are, respectively, the metric projection of $x$ and $y$ onto $\bar{B}(x_0,r)$, then there exists $R > 0$, depending only on $r$ and
$\varepsilon$, such that
$$
\bar{B}(x_0,R)\cap[x,y]\neq\emptyset.
$$
\end{pro}

%%%%%%%%%%%%%%%%%%%%%%%%%%%%%%%%%%%%%%%%%%%%%%%%%%%%%%%%%%%%%%%%%%%%%%%%%%%%%%%%%%%%%

\section{Main result}
\setcounter{equation}{0}\setcounter{thm}{0}

Now we propose our main result.
\begin{thm}
Let $X$ be a~complete CAT(-1) space and a~nonempty $K\subset X$ be closed and convex. Then $K$ has a~fixed point property for nonexpansive mappings $T\colon K\to K$ if and only if $K$ is geodesically bounded.
\end{thm}

{\bf Proof.}

First let us suppose that $K$ is not geodesically bounded, then there exists a~geodesic $c(t)$, $t\geq 0$ and it is easy to construct a~nonexpansive mapping $T\colon K\to l$ which does not have a~fixed point. Indeed, let us define
$$
Tx=c(t+1), \qquad\mbox{when }P_{l}(x)=c(t).
$$

Now let us suppose that $K$ is geodesically bounded and there is a~nonexpansive mapping $T\colon K\to K$ which does not have a~fixed point. We will show that this yields a~contradiction.

\vskip3mm

{\bf Step 1.}

In Step 1 we consider the set of fixed points of contractions based on the mapping $T$.

Let us fixed $\theta\in K$. For each $t\in[0,1)$ we define
$$
T_tx=t\theta+(1-t)Tx.
$$
Clearly, each $T_t$ is a~contraction with $k=(1-t)<1$, so the Banach contraction principle implies that there is a~unique $z_t\in K$ such that $z_t=T_tz_t=t\theta+(1-t)Tz_t$. Next let us choose a~sequence $(t_n)_{n=1}^\infty$ such that $t_n\to 0$ and denote $z_{t_n}$ by $z_n$.

If the sequence $(z_n)$ is bounded then it has a~unique asymptotic center $A((z_n))$ and since $d(z_n,Tz_n) = \dfrac{t_n}{1-t_n}d(\theta,z_n) \to 0$ this asymptotic center is a~fixed point of $T$, a~contradiction. So we obtain $d(\theta,Tz_n) \to \infty$. Without lose a~generality one may suppose that $d(\theta,Tz_n)>n$.

\vskip3mm

{\bf Step 2.}

In Step 2 we consider the behavior of projections of $Tz_n$ onto closed balls $\bar{B}(x_0,m)$, $m\in\mathbb{N}$.

For each $m\in\mathbb{N}$ let us consider a~sequence $(y^m_n)_{n=m}^\infty$ such that $y^m_n\in [\theta,Tz_n]$ and $d(\theta,y^m_n)=m$, $n\geq m$. Now we will show that there is $m\in\mathbb{N}$ such that $(y^m_n)_{n=m}^\infty$ is not totally bounded. Indeed, let us suppose that it is not true. Then we may find a~subsequence of $(Tz_n)$ (denoting again by $(Tz_n)$) such that the sequence $(y^1_n)_{n=1}^\infty$ is a~Cauchy one. Next we find a~subsequence of $(Tz_n)$ (denoting again by $(Tz_n)$) such that $(y^2_n)$ is a~Cauchy sequence. Finally taking a~diagonal sequence we obtain that for each $m\in\mathbb{N}$ the sequence $(y_n^m)$ is a~Cauchy one. Since $K$ is complete as a~closed subset of a~complete space $X$, we have $y_n^m\to y^m\in K$.

Let us consider a~comparison triangles $\Delta(\bar{\theta},\bar{y}^M_n,\bar{y}^M)$ on a~plane of triangles $\Delta(\theta,y^M_n,y^M)$. Then $\bar{y}_n^m$ ($m<M$) -- the comparison point of $y_n^m$ is lying on $[\bar{\theta},\bar{y}_n^M]$. Choosing a~point $u^m\in[\theta,y^M]$ with $d(\theta,u^m)=m$ it is easy to see that $d(y^m_n,u^m)\leq \dfrac{m}{M}d(y^M_n,y^M)$, so $y^m_n\to u^m$ and $y^m\in[\theta,y^M]$ for $m<M$. Since $d(\theta,y^m)=m$ it follows that the limit sequence $(y^m)_{m=1}^\infty$ forms a~geodesic of $X$ (compare the application of totally boundedness of balls in the proof of \cite[Proposition 3.5]{EP}). Since $K$ is convex this geodesic must belong to $K$ and again we obtain a~contradiction with assumptions on $K$. So there is $R\in\mathbb{N}$ such that $(y_n^R)$ is not totally bounded. Moreover, again taking a~subsequence if it is necessary, one may suppose that there is a~positive real number $r$ such that $d(y^R_p,y^R_q)\geq r$, $p,q\geq R$. In the sequel we will consider this
subsequence instead of $(Tz_n)_{n=1}^\infty$.

\vskip3mm

{\bf Step 3.}

Let us notice that since $T$ is nonexpansive and $z_n\in[\theta,Tz_n]$ it follows that
$$
d(\theta,z_n)+d(z_n,Tz_n)=d(\theta,Tz_n)\leq d(\theta,T\theta)+d(T\theta,Tz_n)\leq d(\theta,T\theta)+d(\theta,z_n),
$$
so a~sequence of positive numbers $(d(z_n,Tz_n))$ is bounded. Let us suppose that there is a~subsequence (denoting again by $(Tz_n)$) such that
\begin{equation}\label{dazydod}
d(z_n,Tz_n)\to d>0
\end{equation}
We will show that \eqref{dazydod} leads to a~contradiction.

For each pair $p,q \geq M$ ($p\neq q$ and $p,q$ large enough) we consider $\Delta(\bar{\theta},\bar{T}z_p,\bar{T}z_q)$ -- a~comparison triangle on a~plane of $\Delta(\theta,Tz_p,Tz_q)$.

Let $\bar{u}_p\in [\bar{\theta},\bar{T}z_p]$ ($\bar{u}_q\in [\bar{\theta},\bar{T}z_q]$) be chosen in such a~way that $d(\bar{u}_p,\bar{T}z_p)=d$ ($d(\bar{u}_q,\bar{T}z_q)=d$).
Clearly, comparing angles, it must be
$$
\angle_{\bar{\theta}}(\bar{T}z_p,\bar{T}z_q)=\angle_{\bar{\theta}}(\bar{y}_p^R,\bar{y}_q^R)\geq 2\arcsin \dfrac{r}{2R}.
$$
So on account of Lemma 3.1
$$
d(z_p,z_q)\leq d(\bar{z}_p,\bar{z}_q)\leq d(\bar{u}_p,\bar{z}_p)+d(\bar{u}_p,\bar{u}_q)+d(\bar{z}_q,\bar{u}_q)
$$
$$
\leq d(\bar{u}_p,\bar{z}_p)+d(\bar{z}_q,\bar{u}_q) + d(Tz_p,Tz_q)-d\cdot \dfrac{r^2}{4R^2}.
$$
Since $d(Tz_p,Tz_q) \to \infty$ (compare Lemma 3.1) and $d(\bar{u}_p,\bar{z}_p)$, $d(\bar{z}_q,\bar{u}_q)$ tends to $0$ (if $p,q\to\infty$), for $p,q$ large enough we obtain
$$
d(z_p,z_q)< d(Tz_p,Tz_q),
$$
a~contradiction.
That means that our assumption \eqref{dazydod} cannot hold and $d(z_n,Tz_n)$ must tend to $0$.

\vskip3mm

{\bf Step 4.}

On account of Proposition 3.4 it follows that there is $M>0$ such that each metric segment $[z_p,z_q]$ has a~nonempty intersection with a~closed ball $\bar{B}(\theta,M)$. Let $z_{p,q}$ belong to the intersection $\bar{B}(\theta,M) \cap [z_p,z_q]$ and fix $q=p+1$ for all $p$ large enough ($p\gg M$).
Since $d(\theta,z_n)\to \infty$ when $n\to\infty$, we have that both $d(z_{p,q},z_p)$ and $d(z_{p,q},z_q)$ tend to infinity for $p\to \infty$ (so also $q\to\infty$).
%\begin{figure}[!htb]
%\mbox{\null}\hfill{\psfig{figure=najw3.eps,height=35mm}}\hfill\mbox{\null}
%\caption{Triangle $\Delta(Tz_{p,q},Tz_p,Tz_q)$}
%\end{figure}

Now we estimate the distance $d(z_{p,q},Tz_{p,q})$. To do this let us consider a~triangle $\Delta(Tz_{p,q},Tz_p,Tz_q)$.
Let us denote $d(Tz_n,z_n)=\varepsilon_n$. Then we have
$$
d(Tz_{p,q},Tz_p) \leq d(z_{p,q},z_p),\qquad
d(Tz_{p,q},Tz_q) \leq d(z_{p,q},z_q)
$$
and
$$
d(Tz_p,Tz_q)  \geq d(z_p,z_q)-\big[\varepsilon_p+\varepsilon_q\big].
$$
Let us denote by $u_{p,q}$ a~point of $[Tz_p,Tz_q]$ such that
$$
d(Tz_p,u_{p,q})=d(Tz_p,Tz_q)\cdot\dfrac{d(z_{p,q},z_p)}{d(z_p,z_q)}
$$
and
$$
d(Tz_q,u_{p,q})=d(Tz_p,Tz_q)\cdot\dfrac{d(z_{p,q},z_q)}{d(z_p,z_q)}.
$$

Since
$$
\pi = \angle_{u_{p,q}}(Tz_p,Tz_q)\leq \angle_{u_{p,q}}(Tz_p,Tz_{p,q})+\angle_{u_{p,q}}(Tz_{p,q},Tz_q),
$$
at least one angle of the sum is not smaller than $\pi/2$ (without lose a~generality we assume that in each case that is $\angle_{u_{p,q}}(Tz_p,Tz_{p,q})$).
Hence
$$
d(Tz_p,Tz_{p,q})> d(Tz_p,u_{p,q}) \geq d(z_p,z_{p,q})-\big[\varepsilon_p,\varepsilon_q\big]\cdot \dfrac{d(z_p,z_{p,q})}{d(z_p,z_q)}
$$
$$
\geq d(Tz_p,Tz_{p,q}) -\big[\varepsilon_p,\varepsilon_q\big].
$$
So
$
0< d(Tz_p,Tz_{p,q})-d(Tz_p,u_{p,q})\to 0
$ if $p$ increases.
Since $d(u_{p,q},Tz_p)$ and $d(u_{p,q},Tz_p)$ tend to infinity, on account of Lemma 3.2 we have that
\begin{equation}\label{upqdo0}
d(Tz_{p,q},u_{p,q})\to 0.
\end{equation}

The Busemann convexity of a~CAT(-1) space implies that
$$
d(z_{p,q},u_{p,q})\leq \max\{d(z_p,Tz_p),d(z_q,Tz_q)\},
$$
what on account of \eqref{upqdo0} leads to
\begin{equation}\label{podciag}
d(z_{p,q},Tz_{p,q})\to 0.
\end{equation}
But the sequence $(z_{p,p+1})$ is bounded, so it has a~unique asymptotic center. Using \eqref{podciag} this asymptotic center must be a~fixed point of $T$, what contradicts our assumptions. Let us mention that we suppose that $K$ is geodesically bounded and there is a~nonexpansive mapping $T\colon K\to K$ which does not have a~fixed point. Obtaining a~contradiction means that if $K$ is geodesically bounded each nonexpansive map $T\colon K\to K$ must have a~fixed point, what finishes the proof of our Theorem.$_{}\qquad\square$

\vskip3mm

\begin{rem}
Let us note that we suppose that $\kappa=-1$ only to simplify our estimations. So one may get the same result as above for each space with curvature bounded above by a~negative number $\kappa$ as the following corollary shows.
\end{rem}

\begin{wn}
Let $X$ be a~complete CAT($\kappa$) space with $\kappa<0$ and a~nonempty $K\subset X$ be convex and closed. Then $K$ has a~fixed point property for nonexpansive mappings if and only if $K$ is geodesically bounded.
\end{wn}

Since each $\mathbb{R}$--tree is a~CAT(-1) space, we obtain:

\begin{wn}{\bf (compare \cite{EK} and \cite[Theorem 31]{Kirk})\\}
Let $X$ be a~complete $\mathbb{R}$--tree and a~nonempty $K\subset X$ be convex and closed. Then $K$ has a~fixed point property for nonexpansive mappings if and only if $K$ is geodesically bounded.
\end{wn}

Moreover, in \cite{GR} K.~Goebel and S.~Reich considered the complex Hilbert ball $B$ and the real Hilbert ball being its subset of constant curvature equal to $-1$. Since the sectional curvature of the complex Hilbert space may be estimated by $-1$ and $-4$ (see \cite[Theorem II.10.16]{BH}), we obtain the following generalization of some results proved in \cite{GR}:

\begin{wn}{\bf (compare \cite[Lemma 30.1]{GR} and \cite[Theorem 32.2]{GR})\\}
Let $X$ be the complex Hilbert ball and a~nonempty $K\subset X$ be convex and closed. Then $K$ has a~fixed point property for nonexpansive mappings if and only if $K$ is geodesically bounded.
\end{wn}


\begin{thebibliography}{99}
\bibitem{khamsi}
A.~G.~Aksoy, M.~A.~Khamsi, A~selection theorem in metric trees,
\textit{Proc. Amer. Math. Soc.}, \textbf{134}, 2957--2966 (2006).
\bibitem{BH}
M.~Bridson, A.~Haefliger, Metric spaces of non-positive curvature, Springer-Verlag, Berlin, 1999.
\bibitem{Domingues}
T.~Dom\'inguez Benavides, The failure of the fixed point property for unbounded sets in $c_0$, \textit{ Proc. Amer. Math. Soc.}, \textbf{140}, 645--650 (2012).
\bibitem{EF}
R.~Esp\'inola, A. Fern\'{a}ndez-Le\'{o}n, CAT($\kappa$)-spaces, weak convergence and fixed
points, \textit{J. Math. Anal. Appl.}, \textbf{353}, 410--427 (2009).
\bibitem{EK}
R.~Esp\'inola, W.~A. Kirk, Fixed point theorems in
$\mathbb{R}$-trees with applications to graph theory, \textit{Topology
Appl.}, \textbf{153}, 1046--1055 (2006)
\bibitem{EP}
R.~Esp\'inola, B. Piatek, The fixed point property and unbounded sets in
CAT(0) spaces, \textit{J. Math. Anal. Appl.}, \textbf{408}, 638--654 (2013)
\bibitem{GK}
K.~Goebel, W.~A. Kirk,Topics in metric fixed point theory,  Cambridge studies in advanced mathematics 28, Cambridge University Press, Cambridge, 1990.
\bibitem{GR}
K.~Goebel, S.~Reich, Uniform convexity, hyperbolic geometry and nonexpansive mappings, Pure Appl. Math., Marcel Dekker, Inc., New York--Basel, 1984.
\bibitem{Kirk}
W.~A. Kirk, Geodesic Geometry and Fixed Point Theory, Seminar of Mathematical Analysis (Malaga/Seville, 2002/2003), 195-225, Univ. Sevilla Secr. Publ., Seville, 2003.
\bibitem{mine}
B. Piatek, Best approximation of coincidence points in metric trees,
\textit{ Ann. Univ. Mariae Curie-Sk\l odowska Sect. A}, \textbf{62}, 113--121 (2008).
\bibitem{ORay}
W.O. Ray, The fixed point property and unbounded sets in Hilbert space, \textit{Trans. Amer. Math. Soc.}, \textbf{258}, 531--537 (1980)
\bibitem{Sine}
R. Sine, On the converse of the nonexpansive map fixed point theorem for Hilbert
space, \textit{Proc. Amer. Math. Soc.}, \textbf{100}, 489--490 (1987)
\bibitem{Kohsaka}
W. Takahashi, J.-C. Yao, F.~Kohsaka, The fixed point property and unbounded sets in Banach spaces, \textit{Taiwanese J. Math.}, \textbf{14}, 733--742 (2010)


\end{thebibliography}
\end{document}